\documentclass[12pt]{amsart}
\usepackage{amsfonts,amsmath,amssymb}
\numberwithin{equation}{section}

\newtheorem{theorem}{Theorem}[section]
\newtheorem{lemma}[theorem]{Lemma}
\newtheorem{proposition}[theorem]{Proposition}

\newtheorem{definition}{Definition}[section]

\newtheorem{remark}[theorem]{Remark}

\newcommand{\cl}[1]{\mathcal{#1}} 

\newcommand{\sca}[1]{\left\langle#1\right\rangle} 
\newcommand{\nor}[1]{\left\Vert #1\right\Vert}

\begin{document}

\title{Stably isomorphic dual operator algebras}
\author{G.K. Eleftherakis}
\address{Dept. of Mathematics, University of Athens, Athens, Greece}
\email{gelefth@math.uoa.gr}
\author{V.I. Paulsen}
\address{Dept. of Mathematics, University of Houston, Houston, TX, 77204}
\email{vern@math.uh.edu}
\date{}
\thanks{This project is cofunded by European Social Fund and National
Resources - (EPEAEK II) "Pyhtagoras II" grant No. 70/3/7997}

\maketitle
\begin{abstract} We prove that two unital dual operator algebras $A,B$ are stably isomorphic 
if and only if they are $\Delta $-equivalent \cite{ele2}, if 
and only if they have completely isometric normal representations $\alpha, \beta $ on Hilbert 
spaces $H, K$ respectively and there exists a ternary ring of operators $\cl{M}\subset B(H,K)$ such that 
$\alpha (A)=[\cl{M}^* \beta (B) \cl{M}]^{-w^*}\;\; \text{and}\;\;
\beta (B)=[\cl{M}\alpha (A)\cl{M}^*]^{-w^*}.$

\end{abstract}

\section{Introduction}

Two dual operator algebras \cite{b, lm} $A, B$ are called \textbf{stably isomorphic} 
if there exists a cardinal $I$ such that the algebras $M_I(A), M_I(B)$ of matrices 
indexed by $I,$ whose finite submatrices have uniformly bounded norms, 
are algebraically isomorphic through an isomorphism which is completely 
isometric and $w^*$-(bi)continuous. In the 
special case of $W^*$-algebras \cite{b}, this happens if and only if $A$ and $B$ are Morita 
equivalent in the sense of Rieffel \cite{rif}. A proof of this fact for 
separably acting von Neumann algebras can be found in \cite{ru} 
and the general case is in \cite{b}.

 In \cite{ele1, ele2} two new equivalence relations between dual operator algebras 
were defined:
\begin{definition}\cite{ele1} Let $A, B$ be $w^*$-closed algebras acting on Hilbert spaces 
$H$ and $K,$ respectively. If there exists a \textbf{ternary ring of operators (TRO)} 
$\cl{M}\subset B(H,K),$ i.e. a subspace satisfying $\cl{M}\cl{M}^*\cl{M}\subset \cl{M},$ 
such that 
$A=[\cl{M}^* B \cl{M}]^{-w^*}$ and 
$B=[\cl{M} A \cl{M}^*]^{-w^*}$ we write $A \stackrel{\cl{M}}{\sim} B.$ We say that the algebras 
$A, B$ are \textbf{TRO equivalent} if there exists a TRO $\cl{M}$ such that 
$A \stackrel{\cl{M}}{\sim} B.$ 
\end{definition}

If $A$ is a dual operator algebra, then we call a completely contractive, $w^*$-continuous homomorphism $\alpha : A\rightarrow B(H)$ 
where $H$ is a Hilbert space, a \textbf{normal representation of $A$.}

In \cite{ele2} the notion of $\Delta$-equivalence of two 
unital dual operator algebras $A, B$ was defined in terms of equivalence of two appropriate 
categories. In the present paper, we will adopt the following definition of $\Delta $-

equivalence.

\begin{definition}\label{delta} Two unital dual operator algebras $A, B$ are called 
\textbf{$\Delta$-equivalent} if they have completely isometric normal 
representations $\alpha, \beta $ such that 
the algebras $\alpha(A), \beta(B)$ are TRO equivalent.  
\end{definition}    

\begin{remark}
The conclusion of the present paper (Theorem \ref{32}) was used in  
\cite[Theorem 1.3]{ele2}. It was proved in that theorem that definition \ref{delta} is in fact 
equivalent to the one given in \cite[Definition 1.4]{ele2}: there, 
two unital dual operator algebras $A$ and $B$ are called 
$\Delta$-equivalent if there exists an equivalence functor 
between their categories of normal representations which 
intertwines not only the representations of the algebras 
but also their restrictions to the diagonals. 
\end{remark}

Two completely isometrically and $w^*$-continuously isomorphic unital dual operator algebras 
are not necessarily TRO equivalent, but they are $\Delta$-equivalent. Also two 
$W^*$-algebras are Morita equivalent in the sense of Rieffel if and only if they are 
$\Delta$-equivalent \cite{ele2}. In this work we are going to prove that two unital 
dual operator algebras are $\Delta$-equivalent if and only if they are stably isomorphic.

 We explain now why two stably isomorphic unital dual operator algebras are $\Delta$-equivalent.
We need first to present some definitions and results, see for example \cite{b}. If $I$ is 
a cardinal and $X$ is a dual operator space, we denote by $\Omega _I(X)$ the linear space 
of all matrices with entries in $X.$ If $x\in \Omega _I(X)$ and $r$ is a finite subset of $I$ we 
write $x^r=(x_{ij})_{i,j\in r}.$ We define $$\|x\|=\sup_{r\subset I, finite}\|x^r\|\;\;\text{and}
\;\;M_I(X)=\{x\in \Omega _I(X), \|x\|< +\infty \}.$$
This space is a dual operator space. If $X$ is a dual operator algebra then $M_I(X)$ 
is also a dual operator algebra. In case $X$ is a $w^*$-closed subspace of $B(H,K)$ for some 
Hilbert spaces $H, K$ we naturally identify $M_I(X)$ as a subspace of $B(H^I,K^I)$ where 
$H^I\; (resp. K^I)$ is the direct sum of $I$ copies of $H\; (resp. K).$ We denote the $w^*$-closed 
 subspace of 
$B(H^I,K)$ consisting of bounded operators of the form 
$$H^I\rightarrow K: (\xi _i)_{i\in I}\rightarrow \sum_ix_i(\xi _i)$$ for $\{x_i: i\in I\}\subset 
X$ by $R_I^w(X)$ and the $w^*$-closed subspace of $B(H, K^I)$ consisting of bounded operators of the form 
$$H\rightarrow K^I: \xi \rightarrow (x_i(\xi ))_{i\in I}$$ for $\{x_i: i\in I\}\subset 
X$ by $C_I^w(X).$ Observe that if $X$ is a $w^*$-closed TRO then the spaces $R_I^w(X), 
C_I^w(X)$ are $w^*$-closed TRO's.

Suppose now that the unital dual operator algebras $A_0, B_0$ are stably isomorphic 
for a cardinal $I$. By \cite{lm} there exist completely isometric normal representations of $A_0, B_0$
whose images we denote by $A, B,$ respectively. Observe that the algebras $A, M_I(A)$ are TRO equivalent, indeed,
$A \stackrel{\cl{M}}{\sim} M_I(A),$ where $\cl{M}=C_I^w(\Delta (A)),$ and $\Delta (A)= A \cap A^*$ is the 
diagonal of $A.$ Similarly the algebras $B, M_I(B)$ are TRO equivalent. Since $\Delta $-
equivalence is an equivalence relation preserved by normal completely isometric homomorphisms 
we conclude that the initial algebras are $\Delta $-equivalent.

The purpose of this paper is 
to prove the converse: $\Delta $-equivalent algebras are stably isomorphic. Since 
every completely isometric normal homomorphism $A\rightarrow B$ for dual operator algebras 
naturally \textquotedblleft extends" to a completely isometric normal homomorphism $M_I(A)\rightarrow M_I(B)$ 
for every cardinal $I$ \cite{b}, it suffices to show that the TRO equivalent algebras 
are stably isomorphic.

\section{Generated bimodules.}

In this section we prove that if $A \;(resp. B)$ is a $w^*$-closed subalgebra of $B(H)\; 
(resp. B(K))$ 
for a Hilbert space $H\; (K)$ and $\cl{M}\subset B(H,K)$ is a TRO such that 
$A \stackrel{\cl{M}}{\sim} B,$ then there exist bimodules $X, Y$ over these algebras, i.e., 
$AXB\subset X,\; BYA\subset Y,$ which are generated by $\cl{M},$ such that $A\cong X\stackrel{\sigma h}
{\otimes}_B Y$ and $B\cong Y\stackrel{\sigma h}
{\otimes}_A X$ as dual spaces, where $X\stackrel{\sigma h}
{\otimes}_B Y$ ($Y\stackrel{\sigma h}
{\otimes}_A X$) is an appropriate quotient of the normal Haagerup tensor product $X\stackrel{\sigma h}
{\otimes} Y$ ($ Y\stackrel{\sigma h}
{\otimes} X)$ \cite{er2}.

We start with some definitions and symbols. If $\Omega $ is a Banach space we denote by 
$\Omega ^*$ its dual. 
If $X, Y, Z$ are linear spaces, $n\in \mathbb{N}$ and $\sigma : X\rightarrow Y$ is a linear map 
we denote again 
by $\sigma $ the map $M_n(X)\rightarrow M_n(Y): (x_{ij})\rightarrow (\sigma (x_{ij})).$ 
If $\phi : X\times Y\rightarrow Z$ is a bilinear map and $n, p \in \mathbb{N}$ we denote again 
by $\phi $ the map $M_{n,p}(X)\times M_{p,n}(Y)\rightarrow M_n(Z): ((x_{ij}),(y_{ij}))\rightarrow 
(\sum_{k=1}^p\phi (x_{ik},y_{kj}))_{ij}.$ If $X, Y$ are operator 
spaces we denote by $CB(X,Y)$ the space of completely bounded maps from $X$ to $Y$ 
with the completely bounded norm. If $Z$ is another operator space, a bilinear map 
$\phi : X\times Y\rightarrow Z$ is called completely bounded \cite{paul} if there exists $c>0$ 
such that $\|\phi (x,y)\|\leq c\|x\|\|y\|$ for all $x\in M_{n,p}(X), y\in M_{p,n}(Y), n,p 
\in \mathbb{N}.$ The least such $c$ is the completely bounded norm of $\phi $ 
and it is denoted by $\|\phi \|_{cb}.$ We write $$CB(X\times Y,Z)=\{\phi : 
X\times Y\rightarrow Z,\;\; \phi \;\;\text{is\;\; completely\;\; bounded}\}.$$ 
This is an operator space under the identification $$M_n({CB(X\times Y, Z)})=
CB(X\times Y, M_n(Z))$$ for all $n\in \mathbb{N}.$

We denote the Haagerup tensor product of $X,Y$ by $ X\stackrel{h}{\otimes}Y.$ The map 
$CB(X\times Y,Z)\rightarrow CB(X\stackrel{h}{\otimes}Y,Z): \omega \rightarrow 
\stackrel{\sim }{\omega }$ given by $\stackrel{\sim }{\omega }(x\otimes y)=\omega (x,y)$ 
for all $x\in X, y \in Y$ is a complete isometry.
 If $X,Y$ are dual operator spaces we denote by $CB^\sigma (X,Y)$ the space of 
completely bounded $w^*$-continuous maps. If $Z$ is another dual operator space a bilinear map
 $\phi : X\times Y\rightarrow Z$ is called \textbf{normal} if it is separately $w^*$-continuous. 
We denote by $CB^\sigma (X\times Y,Z)$ the space of completely bounded normal bilinear maps. 

We now recall the normal Haagerup tensor product \cite{er2}. In the rest 
of this section we fix dual operator spaces $X,Y$ and the map 
$$\pi : CB(X\times Y, \mathbb{C})\rightarrow  
CB( X\stackrel{h}{\otimes}Y ,\mathbb{C})=(X\stackrel{h}{\otimes}Y)^*$$ given by  
$\pi (\omega )=\stackrel{\sim }{\omega } , \stackrel{\sim }{\omega }(x\otimes y)=\omega (x,y).$
 We denote by $\Omega _1$ the space $\pi (CB^\sigma (X\times Y,\mathbb{C}))$ and 
by $X\stackrel{\sigma h}{\otimes 
}Y$ the dual of $\Omega _1.$ This space is the $w^*$-closed span of its 
elementary tensors $x\otimes y, x\in X, y\in Y$ and it has the following 
property: For all dual operator spaces $Z$ there exists a complete onto isometry 
$$J:CB^\sigma (X\times Y,Z)\rightarrow CB^\sigma (X\stackrel{\sigma h}{\otimes }Y,Z):
\phi \rightarrow \phi_{\sigma }$$ where $\phi_{\sigma }(x\otimes y)= \phi(x,y).$

We now fix a dual operator algebra $B$ such that $X$ is a right $B$-module and $Y$ 
is left $B$-module and the maps 
$$X\times B\rightarrow X: (x,b)\rightarrow xb,\;\; B\times Y\rightarrow Y: (b,y)\rightarrow by$$ 
are complete contractions and normal bilinear maps. A bilinear map $\omega : X\times Y\rightarrow Z$ 
is called \textbf{$B$-balanced} if $\omega(xb,y)= \omega(x,by) $ for all $x\in X, b\in B, y
\in Y.$ For every dual operator space $Z$ we define the space 
$$ CB^{B\sigma }(X\times Y,Z) =\{\omega \in CB^\sigma (X\times Y,Z): \omega \;\;\text{is\;\;
$B$-balanced}\}.$$
We denote by $\Omega _2$ the space $\pi (CB^{B\sigma }(X\times Y, \mathbb{C})).$ Observe that $\Omega _2$ 
is a closed subspace of $\Omega _1\subset (X\stackrel{h}{\otimes }Y)^*.$ Also we define 
the space $$N=[xb\otimes y-x\otimes by: x\in X, b\in B, y\in Y]^{-w^*}\subset X\stackrel{\sigma h}
{\otimes }Y.$$ We denote by $X\stackrel{\sigma h}{\otimes }_BY$ the space 
$(X\stackrel{\sigma h}{\otimes }Y)/N$ and we use the symbol $x\otimes _By$ for 
$x\otimes y+N, x\in X, y\in Y.$

\begin{proposition}\label{21} The spaces $ X\stackrel{\sigma h}{\otimes }_BY$ and 
 $\Omega _2^*$ are completely 
isometric and $w^*$-homeomorphic. 
\end{proposition}
\textbf{Proof.} The adjoint map $\theta : X\stackrel{\sigma h}{\otimes }Y\rightarrow \Omega _2^* $ 
of the inclusion $\Omega_2\hookrightarrow  \Omega_1 $ is a complete quotient map 
and $w^*$-continuous. Check now that $N=Ker(\theta ). \qquad \Box$

\begin{proposition}\label{22} If $Z$ is a dual operator space and 
$\phi \in CB^{B\sigma }(X\times Y,Z)$ then there exists a $w^*$-continuous and 
completely bounded map $\phi _{B\sigma h}: 
X\stackrel{\sigma h}{\otimes }_BY \rightarrow Z$ such that $\phi _{B\sigma h}(x\otimes _By)=
\phi (x,y)$ for all $x\in X, y\in Y.$  
In fact the map $CB^{B\sigma }(X\times Y,Z)\rightarrow CB^{\sigma }(X\stackrel{\sigma h}
{\otimes }_BY,Z): \phi \rightarrow \phi_{B\sigma h}$ is a complete isometry, onto. 
\end{proposition}  
\textbf{Proof.}  
Suppose that $Z_*$ is the operator space predual of $Z.$ 
 For every $\omega \in Z_*, \omega \circ \phi \in \Omega _2.$ So we can 
define a map $\phi _*: Z_*\rightarrow \Omega _2: \phi _*(\omega )=\omega \circ \phi .$
We denote by $\phi _{B\sigma h}$ the adjoint map of $\phi _*$. So that $\phi_{B\sigma h} \in CB(\Omega_2^*,Z) = CB(X\stackrel{\sigma h}{\otimes}_B Y, Z)$ by Proposition~\ref{21}.  For every $x\in X, y\in Y, \omega \in 
Z_*$ we have $\sca{ \phi _{B\sigma h}(x\otimes _B y) ,\omega }=\sca{\phi (x,y),\omega }$ 
so $\phi _{B\sigma h}(x\otimes _B y) =\phi (x,y).$

Let $i: \Omega_2 \to \Omega_1$ denote the inclusion map so that $q= i^*: \Omega_1^* \to \Omega_2^*$
 is a $w^*$-continuous complete quotient map. The map of composition with $q$ gives a completely 
isometric inclusion, $q^*: CB^{\sigma}(\Omega_2^*, Z) \to CB^{\sigma}(\Omega_1^*, Z).$ 

By Proposition~\ref{21} we may identify $\Omega_2^*= X\stackrel{\sigma h}{\otimes}_B Y$ and also we have $\Omega_1^*= X\stackrel{\sigma h} \otimes Y$ by definition. Thus, modulo these identifications, we have that $q^*: CB^{\sigma}( X\stackrel{\sigma h}{\otimes }_B Y, Z) \to CB^{\sigma}( X\stackrel{\sigma h} \otimes Y, Z)$ is a  $w^*$-continuous complete isometry.

We also have that $CB^{B \sigma}(X \times Y, Z) \subseteq CB^{\sigma}(X \times Y, Z)$ is a 
subspace endowed with the same matrix norms. Thus, $J: CB^{B \sigma}(X \times Y, Z) \to CB^{\sigma}( X\stackrel{\sigma h} \otimes Y, Z)$ is also a completely isometric inclusion.

Now observe that $J(\phi) = q^*(\phi_{B\sigma h}),$ so that $\phi \to \phi_{B \sigma h}$ 
is a complete isometry and $J(CB^{B\sigma}(X \times Y, Z)) \subseteq q^*(CB^{\sigma}( X\stackrel{\sigma h}{\otimes }_B Y, Z)).$  

It remains to show that the map is onto so that the above inclusion is an equality of sets.
To see that $\phi \rightarrow \phi_{B\sigma h}$ is onto 
$CB^{\sigma }(X\stackrel{\sigma h}
{\otimes }_BY,Z),$ let 
$\stackrel{\sim }{\psi }\in CB^{\sigma }( X\stackrel{\sigma h}
{\otimes }_BY ,Z)$ and $\theta: X\stackrel{\sigma h}
{\otimes }Y \rightarrow X\stackrel{\sigma h}
{\otimes }_BY : x\otimes y\rightarrow x\otimes _By $ be the map in Proposition \ref{21}. 
Since $ \stackrel{\sim }{\psi } \circ \theta \in CB^{\sigma }( X\stackrel{\sigma h}
{\otimes }Y ,Z)$ the map $\psi : X\times Y\rightarrow Z$ given by $\psi (x,y)=
\stackrel{\sim }{\psi }\circ \theta(x\otimes y)=\stackrel{\sim }{\psi }(x\otimes _By)$ 
belongs to the space $CB^{\sigma }(X\times Y,Z).$ We have to prove that $\psi $ is 
balanced. If $\omega \in Z_*$ then $\omega \circ \stackrel{\sim }{\psi }$ belongs to 
the predual of $X\stackrel{\sigma h}
{\otimes }_BY.$ So there exists $\chi \in CB^{B\sigma }(X\times Y,\mathbb{C})$ such that 
$\chi (x,y)=\omega (\psi (x,y))$ for all $x\in X, y\in Y.$ Now for every  
$x\in X, y\in Y, b\in B$ we have $$\omega (\psi (xb,y))=\chi (xb,y)=
\chi (x,by)=\omega (\psi (x,by)).$$ The functional $\omega $ is arbitrary in $Z_*$ so 
$\psi (xb,y)=\psi (x,by).$ We have proved that the map $CB^{B\sigma }
(X\times Y,Z)\rightarrow CB^{\sigma } ( X\stackrel{\sigma h}
{\otimes }_BY ,Z): \phi \rightarrow \phi_{B\sigma h}$ is an onto.
$\qquad \Box$

\medskip

Suppose now that $H,K$ are Hilbert spaces, $A$ and $B$ are unital $w^*$-closed subalgebras 
of $B(K)$ and $B(H)$ respectively and $\cl{M}\subset B(K,H)$ is a $w^*$-closed TRO such that 
 $A \stackrel{\cl{M}}{\sim} B.$ 
            
\begin{definition}\label{23} The spaces $ [A\cl{M}^*]^{-w^*} , [\cl{M}A]^{-w^*} $ 
are called the \textbf{$\cl{M}$-generated $A-B$ bimodules.} 
\end{definition}

In what follows we assume that $ X=[A\cl{M}^*]^{-w^*}, Y=[\cl{M}A]^{-w^*}. $ We can check that 
$$X=[\cl{M}^*B]^{-w^*}, Y=[B\cl{M}]^{-w^*},$$
\begin{equation}\label{sxeseis}AXB\subset X,\; BYA\subset Y,\; A=[XY]^{-w^*},\; B=[YX]^{-w^*}.
\end{equation}

Let $a\in A.$ We define a map $$ CB^{B\sigma }(X\times Y,\mathbb{C}) \rightarrow 
CB^{B\sigma }(X\times Y,\mathbb{C}): \omega \rightarrow \omega _a, $$
by $\omega _a(x,y)=\omega (x,ya).$ This map is continuous. The adjoint map 
$\pi _a: X\stackrel{\sigma h}{\otimes}_B Y \rightarrow X\stackrel{\sigma h}{\otimes}_B Y $ 
satisfies $\pi _a(x\otimes_ By)=x\otimes _B(ya).$ For every $z\in 
X\stackrel{\sigma h}{\otimes}_B Y $ we define $za=\pi _a(z).$ Observe that if 
$\left(\sum_{i=1}^{k_j}x_i^{j}\otimes _By_i^j\right)_j$ is a net such that $z=w^*-\lim_j\sum_{i=1}^{k_j}
x_i^j\otimes 
_By_i^j$ then $za=w^*-lim_j\sum_{i=1}^{k_j}x_i^j\otimes 
_B(y_i^ja).$ 

\begin{lemma}Let $z\in X\stackrel{\sigma h}{\otimes}_B Y .$ If $(a_\lambda )_\lambda \subset A$ 
is a net such that $a_\lambda \stackrel{w^*}{\rightarrow }a$ then 
$za_\lambda \stackrel{w^*}{\rightarrow }za.$ 
\end{lemma} 
\textbf{Proof.} Choose $\omega \in Ball(CB^\sigma (X\times Y, \mathbb{C})).$ From the 
 normal version of the Christensen, Sinclair, Paulsen, Smith theorem, 
see for example Theorem 5.1 in \cite{er2}, there exist a 
Hilbert space $H$ and normal completely contractive maps $\phi _1: X\rightarrow B(H,\mathbb{C}),$    
$\phi _2: Y\rightarrow B(\mathbb{C},H)$ such that $\omega (x,y)=\phi_1(x)\phi_2(y).$ Observe that 
the bilinear map $Y\times A\rightarrow B(\mathbb{C},H): (y,a)\rightarrow \phi _2(ya)$ is 
completely contractive and normal. So by the same theorem there exist a Hilbert 
space $K$ and complete contractions $\phi _3: A\rightarrow B(\mathbb{C},K), \phi_4:
 Y\rightarrow B(K,H)$ such that $\phi _2(ya)=\phi_4(y) \phi_3(a)$ for all 
$y\in Y, a\in A.$ The bilinear map $X\times Y\rightarrow B(K,\mathbb{C}): (x,y)\rightarrow 
\phi_1(x) \phi_4(y)$ is normal and a complete contraction. So there exists a completely 
contractive $w^*$-continuous map $\pi : X\stackrel{\sigma h}{\otimes}Y \rightarrow B(K,\mathbb{C})$ 
such that $\pi (x\otimes y)=\phi_1(x) \phi_4(y).$ Now the map 
$$\tau (\omega ): (X\stackrel{\sigma h}{\otimes}Y) \times A\rightarrow \mathbb{C}: \tau (\omega )
(z,a)=\pi (z)\phi _3(a)$$ is normal, completely contractive and satisfies 
\begin{align*}&\tau (\omega )(x\otimes y, a)=\pi (x\otimes y)\phi _3(a)\\
=&\phi_1(x) \phi_4(y) \phi_3(a)=\phi_1(x) \phi_2(ya)=\omega (x,ya)  
\end{align*}
for all $x\in X, y\in Y, a\in A.$ The conclusion is that we can define a contraction 
$$\tau : CB^\sigma (X\times Y,\mathbb{C})\rightarrow CB^\sigma (X\stackrel{\sigma h}{\otimes}
Y\times A, \mathbb{C} ): \omega \rightarrow \tau (\omega) $$
which has adjoint map $\sigma : (X\stackrel{\sigma h}{\otimes}Y)\stackrel{\sigma h}{\otimes }A
\rightarrow X\stackrel{\sigma h}{\otimes}Y$ satisfying $\sigma ((x\otimes y)\otimes a)=
x\otimes (ya).$ We recall from Proposition \ref{21} the map 
$$\theta :X\stackrel{\sigma h}{\otimes}Y\rightarrow X\stackrel{\sigma h}{\otimes}_B Y: 
\theta (x\otimes y)=x\otimes _By.$$ 
Choose arbitrary $z\in X\stackrel{\sigma h}{\otimes}_B Y $ and $z_0\in 
X\stackrel{\sigma h}{\otimes}Y $ such that $\theta (z_0)=z.$ If 
$\left(\sum_{i=1}^{k_j}x_i^j\otimes y_i^j\right)_j$ is a net such that $z_0=w^*- lim\sum_{i=1}^{k_j}
x_i^j\otimes 
y_i^j $ then for all $a\in A$
\begin{align*}\theta \circ \sigma (z_0\otimes a)=&\theta \circ \sigma 
\left(lim_j\left( \left(\sum_{i=1}^{k_j}x_i^j\otimes 
y_i^j\right)\otimes a\right) \right)=lim_j\sum_{i=1}^{k_j}\theta (x_i^j\otimes 
(y_i^ja))\\=&lim_j\sum_{i=1}^{k_j}x_i^j\otimes_B 
(y_i^ja)=za.  
\end{align*}         
If $(a_\lambda )_\lambda \subset A$ 
is a net such that $a_\lambda \stackrel{w^*}{\rightarrow }a$ then 
$z_0\otimes a_\lambda \stackrel{w^*}{\rightarrow }z_0\otimes a$ in 
$(X\stackrel{\sigma h}{\otimes}Y)\stackrel{\sigma h}{\otimes }A.$ Since 
$\theta \circ \sigma $ is $w^*$-continuous we have 
$\theta \circ \sigma (z_0\otimes a_\lambda) \stackrel{w^*}{\rightarrow }\theta \circ \sigma 
(z_0\otimes a)$ or equivalently $za_\lambda \stackrel{w^*}{\rightarrow }za. \qquad \Box$

\begin{theorem}\label{24}$A\cong X\stackrel{\sigma h}
{\otimes}_B Y $ and $B\cong Y\stackrel{\sigma h}
{\otimes}_A X$ completely isometrically and $w^*$-homeomorphically.
\end{theorem} 
\textbf{Proof.} The map $X\times Y\rightarrow A: (x,y)\rightarrow xy$ is normal, 
completely contractive and $B$-balanced. So by Proposition \ref{22} it defines a completely 
contractive and $w^*$-continuous map $$\pi : X\stackrel{\sigma h}
{\otimes}_B Y \rightarrow A: \pi (x\otimes _By)=xy.$$ We shall show that $\pi $ 
is a complete isometry. Since $A=[XY]^{-w^*},$ it will follow from the Krein Smulian theorem 
that $\pi $ is onto $A.$ 

Let $z=(z_{ij})\in M_n(X\stackrel{\sigma h}
{\otimes}_B Y ).$ It suffices to show that $\|z\|\leq \|\pi (z)\|.$ Since 
$X\stackrel{\sigma h}
{\otimes}_B Y =( CB^{B\sigma }(X\times Y, \mathbb{C}) )^*$ given $\epsilon >0$ 
 there exist $m\in \mathbb{N}$ 
and $(\omega _{kl})\in Ball(M_m(CB^{B\sigma }(X\times Y, \mathbb{C}) ))$ such that 
$$\|z\|-\epsilon <\|((\omega _{kl}(z_{ij}))_{ij})_{kl}\|.$$ 
By Lemma 8.5.23 in \cite{b} there exist partial isometries $\{v_i: i\in I\}\subset \cl{M}$ 
with mutually orthogonal initial spaces such that $I_H=\sum_{i\in I}\oplus v_i^*v_i.$ By 
the above lemma $$w^*-\lim_{{F\subset I}\atop{finite}}\sum_{s\in F}z_{ij}v_s^*v_s=z_{ij}$$
so $$\lim_{{F\subset I}\atop{finite}}\sum_{s\in F}\omega _{kl}(z_{ij}v_s^*v_s)=\omega _{kl}(z_{ij})
$$ for all $k,l,i,j.$ It follows that there exist partial isometries $\{v_1,...,v_r\}\subset 
\cl{M}$ such that $$\|z\|-\epsilon \leq \|((\sum_{s=1}^r\omega _{kl}(z_{ij}v_s^*v_s
))_{ij})_{kl}\|.$$ Since $X\stackrel{\sigma h}
{\otimes}_B Y$ is the $w^*$-closure of the space $(X\otimes Y)/N,$ see Proposition \ref{21}, 
there exists a net $(z_\lambda )_\lambda \subset M_n(X\otimes Y/N)$ such that $z_\lambda 
\stackrel{w^*}{\rightarrow }z.$ If $z_\lambda =(z_{ij}(\lambda ))_{ij}$ for all $\lambda $ 
we have $z_{ij}(\lambda) 
\stackrel{w^*}{\rightarrow }z_{ij},$ hence $ \sum_{s=1}^r\omega _{kl}
(z_{ij}(\lambda )v_s^*v_s) \rightarrow \sum_{s=1}^r\omega _{kl}
(z_{ij}v_s^*v_s) $ for all $i,j,k,l.$ It follows that there exists $\lambda _0$ such that 
$$\|z\|-\epsilon \leq \nor{ \left(\left( \sum_{s=1}^r \omega _{kl}(z_{ij}(\lambda )v_s^*v_s
) \right)_{ij}\right)_{kl}} \;\;\text{for\;\; all}\;\;\lambda \geq \lambda_0.$$ Fix $i,j,\lambda $ and suppose that 
$z_{ij}(\lambda )=\sum_{p=1}^tx_p\otimes _By_p,$ then  
$\omega _{kl}(z_{ij}(\lambda )v_s^*v_s)=\sum_{p=1}^t\omega _{kl}(x_p,y_pv_s^*v_s)$ for all 
$k,l,s.$ Since $y_pv_s^*\in YX\subset B$ and $\omega _{kl}$ is $B$-balanced we have 
$$ \omega _{kl}(z_{ij}(\lambda )v_s^*v_s) =\sum_{p=1}^t\omega _{kl}(x_py_pv_s^*,v_s)=
\omega _{kl}(\pi (z_{ij}(\lambda ))v_s^*,v_s).$$ 
So we take the inequality $$\|z\|-\epsilon \leq \nor{\left(\left( \sum_{s=1}^r \omega _{kl}
(\pi (z_{ij}(\lambda) )v_s^*,v_s
) \right)_{ij}\right)_{kl}}\;\;\text{for\;\; all}\;\;\lambda \geq \lambda_0.$$ 
Since $\pi(z_{ij}(\lambda ))\stackrel{w^*}
{\rightarrow } \pi(z_{ij})$ we have $$\|z\|-\epsilon \leq \nor{\left(\left( \sum_{s=1}^r \omega _{kl}
(\pi (z_{ij})v_s^*,v_s) \right)_{ij}\right)_{kl}}_{mn}.$$
Let $v=(v_1,...,v_r)^t$ and $$x=(\pi (z_{ij}))_{ij}\cdot
\left[\begin{array}{clr}
v^* &  &  \\ 0 & \ddots &  0 \\ & & v^* \end{array}\right]
\in M_{n,nr}(X),\;\;y=\left[\begin{array}{clr}
v &  &  \\ 0 & \ddots &  0 \\ & & v \end{array}\right]\in M_{nr,n}(Y).$$ 
The above inequality can be written in the following form
$$\|z\|-\epsilon \leq \|(\omega _{kl}(x,y))_{k,l}\|_{mn}.$$ Since 
$$\|(\omega _{kl})\|_m=\|(\omega _{kl}): X\times Y\rightarrow M_m\|_{cb}\leq 1$$ we have 
$$\|z\|-\epsilon \leq \|x\|\|y\|\leq \|(\pi (z_{ij}))_{ij}\|\|v^*\|\|v\|\leq \|\pi (z)\|.$$ 
Since $\epsilon > 0 $ is arbitrary we obtain $\|z\|\leq \|\pi (z)\|.$ This completes the proof of 
$A\cong X\stackrel{\sigma h}
{\otimes}_B Y .$ Similarly we can prove $B\cong Y\stackrel{\sigma h}
{\otimes}_A X \qquad \Box$

\section{The main theorem}  

In this section we shall prove that two unital dual operator algebras are $\Delta $-equivalent 
if and only if they are stably isomorphic. As we noted in section 1 it suffices to show 
that TRO equivalent algebras are stably isomorphic. Thus in what follows, we fix unital 
$w^*$-closed algebras $A,B$ acting on Hilbert spaces $H,K$ respectively and a $w^*$-closed TRO 
$\cl{M}$ such that $A\stackrel{\cl{M}}{\sim }B.$ Let $X=[A\cl{M}^*]^{-w^*} , Y=[\cl{M}A]^{-w^*}$ 
 be the $\cl{M}$-generated $A-B$ 
bimodules which satisfy (\ref{sxeseis}). We 
give the following definition (see the analogous definition in \cite{bmp}). If 
$U_i\subset B(L,H), V_i\subset B(H,L), i=1,2$ are spaces such that $U_iV_i\subset A, i=1,2$ 
a pair of maps $\sigma : U_1\rightarrow U_2,  \pi : V_1\rightarrow V_2$ is called 
\textbf{A-inner product preserving} if $\sigma(x) \pi(y)=xy$ for all $x\in U_1, 
y\in V_1.$

\begin{lemma}\label{31} There exist a cardinal $I$ and completely isometric, $w^*$-
continuous, onto, $A$-module maps $\sigma: R_I^w(X)\rightarrow R_I^w(A), \pi: C_I^w(Y)\rightarrow 
C_I^w(A) $ such that the pair $(\sigma, \pi )$ is $A$-inner product preserving. 
\end{lemma}
\textbf{Proof.} From Lemma 8.5.23 in \cite{b} there exist partial isometries 
$\{m_i: i\in I \}\subset \cl{M}$ with mutually orthogonal initial spaces and 
$\{n_j: j\in J \}\subset \cl{M}$ with mutually orthogonal final spaces such that 
$\sum_{i\in I}\oplus m^*_im_i=I_H,\,\sum_{j\in J}\oplus n^*_jn_j=I_K.$

By introducing sufficiently many 0 partial isometries to each set, we may assume 
that $I^2=I=J.$ We denote by $m$ the column $(m_i)_{i\in I}\in C_I^w(\cl{M}).$ 
We have $m^*m=I_H$ and we denote by $p$ the projection $mm^*\in M_I(B).$

In what follows if $U_n\subset B(H_n,K)$ are $w^*$-closed subspaces, $H_n, K$ 
Hilbert spaces, $n\in \mathbb{N},$ we denote by $U_1\oplus _rU_2\oplus _r...$ the 
$w^*$-closed subspace of $B(\sum_n\oplus H_n,K)$ generated by the bounded 
operators of the form $(u_1,u_2,...),$$ u_n\in U_n, n\in \mathbb{N}.$ Also 
if $V_n\subset B(K,H_n)$ are $w^*$-closed subspaces, $H_n, K$ 
Hilbert spaces, $n\in \mathbb{N}$ we denote by $V_1\oplus _cV_2\oplus _c...$ the 
$w^*$-closed subspace of $B(K,\sum_n\oplus H_n)$ generated by the bounded 
operators of the form $(v_1,v_2,...)^t,$$ v_n\in V_n, n\in \mathbb{N}.$
If $(x_i)_{i\in I}\in R_I^w(R_I^w(X)) $ where $x_i\in R_I^w(X)$ then $x_im\in A$ 
and so we can define the maps 
$$\tau _1: R_I^w(R_I^w(X)) \rightarrow R_I^w(A)\oplus _rR_I^w(R_I^w(X)p^\bot 
), $$$$\tau _1( (x_i)_{i\in I} )=((x_im)_{i\in I},(x_ip^\bot )_{i\in I}),\;\; x_i \in R_I^w(X) 
$$ and 
$$\tau _2: C_I^w(C_I^w(Y)) \rightarrow C_I^w(A)\oplus _c
C_I^w(p^\bot C_I^w(Y)), $$$$\tau _2( (y_i)_{i\in I})=((m^*y_i)_{i\in I},(p^\bot y_i )_{i\in I})^t, 
\;\;y_i \in C_I^w(Y).$$
This pair of maps is $A$-inner product preserving: if $x\in R_I^w(R_I^w(X)), y\in C_I^w(C_I^w(Y))$ 
then $$\tau_1(x)\tau_2(y)=(xm,xp^\bot )(m^*y,p^\bot y)^t=xmm^*y+xp^\bot y=xpy+xp^\bot y=xy.$$
These maps are onto because every $a\in A$ may be written $a=(am^*)m$ with $am^*\in R_I^w(X)$  
 and also $a=m^*(ma)$ with $ma\in C_I^w(Y)$ and they are clearly $w^*$-continuous $A$-module maps. 
Also they are complete isometries. We check this fact for $\tau _1$ and $n=2:$ If 
$x=(x_{ij})\in M_2(R_I^w(R_I^w(X)))$ we have 
\begin{align*} &\|\tau _1(x)\|^2= \nor{\left[ \begin{array}{cllr} x_{11}m & 
x_{11}p^\bot & x_{12}m & x_{12}p^\bot \\ x_{21}m & x_{21}p^\bot & x_{22}m & 
x_{22}p^\bot \end{array}\right]}^2 \\=&\nor{\left[\begin{array}{cllr} x_{11}m & 
x_{12}m & x_{11}p^\bot  & x_{12}p^\bot \\ x_{21}m & x_{22}m & x_{21}p^\bot  & 
x_{22}p^\bot \end{array}\right]}^2
= \nor{\left[ x \left[\begin{array}{clr}m & 0\\ 0 & m \end{array}\right] , 
x \left[\begin{array}{clr}p^\bot  & 0\\ 0 & p^\bot  \end{array}\right] \right]}^2\\=&
\nor{ x \left[ \begin{array}{clr}m & 0\\ 0 & m \end{array}\right] 
\left[ \begin{array}{clr}m^* & 0\\ 0 & m^* \end{array}\right]   
x^* + x\left[ \begin{array}{clr}p^\bot  & 0\\ 0 & p^\bot \end{array}\right] x^*}^2=\|xx^*\|=
\|x\|^2.
\end{align*}
We use the symbol $\infty $ for the $\aleph_0$ cardinal. The following spaces are isomorphic 
as $A$-modules and as dual operator spaces:
\begin{align*}R_\infty ^w(R_I^w(R_I^w(X)))\cong & R_I^w(A)\oplus _rR_I^w(R_I^w(X)p^\bot )\oplus _r
R_I^w(A)\oplus _r...\\\cong &R_I^w(A)\oplus _rR_\infty ^w(R_I^w(R_I^w(X)))
\end{align*}  
 
and \begin{align*}C_\infty ^w(C_I^w(C_I^w(Y)))\cong & C_I^w(A)\oplus _cC_I^w(p^\bot C_I^w(YX)
)\oplus _cC_I^w(A)\oplus _c...\\\cong &C_I^w(A)\oplus _cC_\infty ^w(C_I^w(C_I^w(Y)))
\end{align*}

Since $I^2=I$ it follows that $\infty I=I$ so we have $$ R_I^w(X)\cong R_\infty ^w(R_I^w(R_I^w(X))) \;\;
\text{and}\;\; C_I^w(Y)\cong C_\infty ^w(C_I^w(C_I^w(Y))).$$
 We conclude that there exist completely isometric, $w^*$-continuous, $A$-module bijections 
$$ \lambda _1: R_I^w(X)\rightarrow R_I^w(A)\oplus _rR_I^w(X) \;\;\text{and}
\;\;\lambda _2: C_I^w(Y)\rightarrow C_I^w(A)\oplus _cC_I^w(Y) .$$
 We can choose $\lambda_1, \lambda_2 $ to be $A$-inner product preserving. Similarly 
working with the partial isometries $\{n_j: j\in I\}$ (see the beginning of the proof)
 we obtain a pair 
$(\nu _1,\nu _2)$ of $A$-inner product preserving, completely isometric, $w^*$-continuous 
$A$-module bijections: 
$$ \nu _1: R_I^w(A)\oplus _rR_I^w(X)\rightarrow R_I^w(A) \;\;\text{and}
\;\;\nu _2: C_I^w(A)\oplus _cC_I^w(Y)\rightarrow C_I^w(A) .$$
The maps $$\sigma =\nu_1\circ  \lambda_1: R_I^w(X)\rightarrow R_I^w(A)\;\;\text{and}
\;\;\pi= \nu_2\circ  \lambda_2: C_I^w(Y)\rightarrow C_I^w(A)$$ 
satisfy our requirements. $\qquad \Box$      
  
\begin{theorem}\label{32} Two unital dual operator algebras are $\Delta -equivalent$ 
if and only if they are stably isomorphic.
\end{theorem}
\textbf{Proof.} It suffices to show that if the algebras, $A$ and $B,$ are TRO-equivalent, then they are stably isomorphic. Let 
$I,\sigma, \pi $ be as in Lemma \ref{31}. Observe that $A\stackrel{C_I^w(\cl{M})}{\sim }M_I(B)$ 
and the $C_I^w(\cl{M})$-generated $A-M_I(B)$ bimodules (see definition \ref{23}) 
are the spaces $R_I^w(X)$ and $C_I^w(Y).$ So by Theorem \ref{24} the map $$\psi _1: 
C_I^w(Y)\stackrel{\sigma h}{\otimes}_AR_I^w(X)\rightarrow M_I(B): \psi _1(y\otimes _Ax)=yx$$ 
is a completely isometric, $w^*$-continuous bijection. For the same reason 
the map $$\psi _2: 
C_I^w(A)\stackrel{\sigma h}{\otimes}_AR_I^w(A) \rightarrow M_I(A): \psi _2(a\otimes _Ac)=ac$$ 
is a completely isometric, $w^*$-continuous bijection. The map $$C_I^w(Y)\times R_I^w(X)\rightarrow 
C_I^w(A)\stackrel{\sigma h}{\otimes}_AR_I^w(A): (y,x)\rightarrow \pi(y)\otimes _A \sigma(x)$$ 
is completely contractive, separately $w^*$-continuous and $A$-balanced. So by Proposition 
\ref{22} there exists a completely contractive $w^*$-continuous map 
$$ C_I^w(Y)\stackrel{\sigma h}{\otimes }_AR_I^w(X) \rightarrow 
C_I^w(A)\stackrel{\sigma h}{\otimes}_AR_I^w(A) : y\otimes _Ax\rightarrow 
\pi(y)\otimes _A \sigma(x).$$ 
We denote this map by $\pi \otimes \sigma.$ Similarly we can define 
a complete contraction $\pi^{-1} \otimes \sigma^{-1}: 
C_I^w(A)\stackrel{\sigma h}{\otimes}_AR_I^w(A)\rightarrow  
C_I^w(Y)\stackrel{\sigma h}{\otimes }_AR_I^w(X).$ Since $\pi^{-1} \otimes \sigma^{-1}$ 
is the inverse of $\pi \otimes \sigma $ we conclude that $\pi \otimes \sigma $ 
is a complete isometry. It follows that the map 
$$\gamma =\psi _2\circ (\pi \otimes \sigma )\circ \psi_1^{-1}: M_I(B)\rightarrow M_I(A)$$
 is a completely isometric, $w^*$-continuous bijection. 
It remains to check that it is an algebraic homomorphism. 
Since $M_I(B)=[C_I^w(Y)R_I^w(X)]^{-w^*}$ it suffices to show that $ \gamma(y_1x_1\cdot y_2x_2) =
 \gamma(y_1x_1)\cdot \gamma(y_2x_2)$ for all $x_1, x_2 \in R_I^w(X), y_1, y_2\in C_I^w(Y).$
Indeed, 
\begin{align*}&\gamma(y_1x_1 \cdot y_2x_2)\\= &\psi_2\circ (\pi \otimes \sigma )\circ \psi _1^{-1}
(y_1x_1y_2\cdot x_2)=\;\;\;\; (y_1x_1y_2\in C_I^w(Y), x\in R_I^w(X)) 
\\=&\psi _2\circ (\pi \otimes \sigma )(y_1x_1y_2\otimes _Ax_2)=\psi _2(\pi (y_1x_1y_2)\otimes _A
\sigma (x_2))\\
=&\pi (y_1x_1y_2)\sigma (x_2)=\qquad  (x_1y_2\in A\;\;\text{ and}\;\; \pi\;\;\text{ is\;\;a\;\; 
A-module\;\; map} )
\\=&\pi (y_1)x_1y_2\sigma (x_2)= \qquad((\sigma, \pi ) \;\;\text{is\;\; A-inner\;\; product\;\;
 preserving})\\
=&\pi(y_1) \sigma(x_1) \pi(y_2) \sigma(x_2)= \psi_2(\pi(y_1)\otimes _A \sigma(x_1) )\cdot
\psi_2(\pi(y_2)\otimes _A \sigma(x_2) )\\
=& 
\psi _2\circ (\pi \otimes \sigma )(y_1\otimes _Ax_1)\cdot
\psi _2\circ (\pi \otimes \sigma )(y_2\otimes _Ax_2)
\\
=&\psi _2\circ (\pi \otimes \sigma )\circ \psi ^{-1}(y_1x_1)\cdot
\psi _2\circ (\pi \otimes \sigma )\circ \psi ^{-1}(y_2x_2)=\gamma(y_1x_1)\cdot \gamma(y_2x_2)                    
\end{align*}$\qquad \Box$

\begin{remark} When the  unital dual operator algebras $A,B$ have
  completely isometric normal representations $\alpha, \beta$ on
  separable, Hilbert spaces such that $\alpha(A)$ and $\beta(B)$ are
  TRO equivalent, then the proof of the above theorem shows that
  $M_{\infty}(A)$ and $M_{\infty}(B)$ are completely isometrically isomorphic, i.e., the index
  set $I$ may be taken to be countable.
\end{remark}

\section{Stably isomorphic CSL algebras.}

In this section we assume that all Hilbert spaces are separable. A set of projections on a 
Hilbert space is called a \textbf{lattice} if it contains the zero and identity operators 
and is closed under arbitrary suprema and infima. If $A$ is a subalgebra of $B(H)$ for some Hilbert 
space $H,$ the set 
$$\mathrm{Lat}(A)=\{l\in pr(B(H)): l^\bot Al=0\}$$ is a lattice. Dually if $\cl{L}$ is a lattice 
the space $$\mathrm{Alg}(\cl{L}) =\{a\in B(H): l^\bot al=0 \;\;\forall\;\; l\in \cl{L}\}$$ 
is an algebra. A commutative subspace lattice \textbf{(CSL)} is a projection lattice $\cl{L}$ 
whose elements commute; the algebra  $\mathrm{Alg}(\cl{L})$ is called a \textbf{CSL algebra}. 

Let $\cl{L}$ be a CSL and $l\in \cl{L}.$ We denote by $l_\flat $ the projection 
$\vee \{r\in \cl{L}: r<l\}.$ Whenever $l_\flat <l$ we call the projection $l-l_\flat $ 
an \textbf{atom} of $\cl{L}.$ If the CSL $\cl{L}$ has no atoms we say that it is a \textbf{
continuous CSL}. If the atoms span the identity operator we say that $\cl{L}$ is a 
\textbf{totally atomic CSL}.

If $\cl{L}_1, \cl{L}_2$ are CSL's, $\phi : \cl{L}_1\rightarrow \cl{L}_2$ is a \textbf{lattice 
isomorphism} (a bijection which preserves order) and $p$ (resp. $q$) is the span 
of the atoms of $\cl{L}_1$ (resp. of $\cl{L}_2$) there exists a well defined lattice isomorphism 
$\cl{L}_1|_p\rightarrow \cl{L}_2|_q: l|_p\rightarrow \phi (l)|_q$ (Lemma 5.3 in \cite{ele1}.)
 Observe that the CSL's $\cl{L}_1|_{p^\bot },\;\; $ $\cl{L}_2|_{q^\bot }$ are continuous. 
But it is not always true that $\phi $ induces a lattice isomorphism from $\cl{L}_1|_{p^\bot }$ 
onto $\cl{L}_1|_{q^\bot }.$ In \cite[7.19]{dav} there exists an example 
of isomorphic nests $\cl{L}_1, \cl{L}_2$ such that $p^\bot =0$ and $q^\bot \neq 0.$ This 
motivates the following definition: 

\begin{definition}\label{5.2.d}\cite{ele1} Let $\cl{L}_1, \cl{L}_2$ be CSL's, 
$\phi : \cl{L}_1\rightarrow \cl{L}_2$ be a lattice isomorphism, 
$p$ the span of the atoms of $\cl{L}_1$ and $q$ the span of the atoms of 
$\cl{L}_2.$ We say that $\phi $ \textbf{respects continuity} if there 
exists a lattice isomorphism $\cl{L}_1|_{p^\bot} \rightarrow \cl{L}_2|_{q^\bot}$ 
such that $l|_{p^\bot} \rightarrow \phi (l)|_{q^\bot} $ for every $l\in \cl{L}_1.$
\end{definition}

The following was proved in \cite{ele1} (Theorem 5.7).

\begin{theorem}\label{th1} Let $\cl{L}_1, \cl{L}_2$ be separably acting CSL's. The algebras 
$\mathrm{Alg}(\cl{L}_1), \;$$\mathrm{Alg}(\cl{L}_2)$ are TRO equivalent if and only if 
there exists a lattice isomorphism $\phi : \cl{L}_1\rightarrow \cl{L}_2$ which respects continuity.  
\end{theorem}

Also we recall Theorem 3.2 in \cite{ele3}.

\begin{theorem}\label{th2} Two CSL algebras are $\Delta $-equivalent if and only if they are TRO equivalent.
\end{theorem}

Combining Theorems \ref{th1}, \ref{th2} with Theorem \ref{32} we obtain the following:

\begin{theorem}Two CSL algebras, acting on separable Hilbert 
spaces, are stably isomorphic if and only if there exists a lattice isomorphism 
between their lattices which respects continuity. 
\end{theorem}          

\begin{remark} In fact, since the CSL algebras, say $Alg(\cl{L}_i), i=1,2$ are acting on separable
  Hilbert spaces, we have that if there exists a lattice isomorphism
  between $\cl{L}_1$ and $\cl{L}_2$ that respects continuity, then
  $M_{\infty}(Alg(\cl{L}_1))$ and $M_{\infty}(Alg(\cl{L}_2))$ are
  completely isometrically isomorphic.
\end{remark}   
A consequence of this theorem is that two separably acting CSL algebras with continuous 
or totally atomic lattices are stably isomorphic if and only if they have isomorphic lattices.

\end{document}